\documentclass[12pt]{amsart}
\usepackage{amsmath}
\usepackage{amssymb}
\usepackage{epsfig}
\newcommand{\ignore}[1]{\relax}

\newcommand{\R}{\mathbb R}
\newcommand{\Z}{\mathbb Z}

\newcommand{\M}{\mathcal M}

\newcommand{\ft}{\operatorname{ft}}

\newcommand{\OO}{\mathcal{O}}

\newtheorem{lem}{Lemma}[section]

\newtheorem{cor}[lem]{Corollary}

\newtheorem{theorem}{Theorem}

\newtheorem{prop}[lem]{Proposition}

\theoremstyle{definition}
\newtheorem{defn}[lem]{Definition}
\newtheorem{property}[lem]{Property}

\theoremstyle{remark}
\newtheorem{rmk}[lem]{Remark}

\newtheorem{ack}{\bf Acknowledgements}

\newcommand{\dd}{\partial}

\newcommand{\T}{\mathbb{T}}

\renewcommand{\setminus}{\smallsetminus}

\newcommand{\bMfive}{\overline{\mathcal M}_{0,5}}
\newcommand{\MM}{{\mathcal M}}

\begin{document}
\title{Moduli spaces of rational tropical curves}
\author{Grigory Mikhalkin}
\address{Department of Mathematics, University of Toronto,
40 St George St, Toronto ON M5S 2E4 Canada}

\begin{abstract}
This note is devoted to the definition
of moduli spaces of rational tropical curves with
$n$ marked points. We show that this space has
a structure of a smooth tropical variety of dimension
$n-3$. We define the Deligne-Mumford compactification
of this space and tropical $\psi$-class divisors.
\end{abstract}
\maketitle

This paper gives a detailed description of the moduli
space of tropical rational curves mentioned in \cite{Mi-ICM}.
The survey \cite{Mi-ICM} was prepared
under rather sharp time and volume constraints.
As a result the coordinate presentation of this moduli
space from \cite{Mi-ICM} contains a mistake (it was over-simplified).
In this paper we'll correct the mistake and give
a detailed description on $\bMfive$ as our main example.

\section{Introduction: smooth tropical varieties}
In this section we follow the definitions of \cite{Mi-book}
and \cite{Mi-ICM}.

The underlying algebra of tropical geometry is given
by the semifield $\T=\R\cup\{-\infty\}$ of tropical numbers.
The tropical arithmetic operations are $``a+b"=\max\{a,b\}$ and
$``ab"=a+b$. The quotation marks are used to distinguish between
the tropical and classical operations.
With respect to addition $\T$ is a commutative semigroup with zero $``0_{\T}"=-\infty$.
With respect to multiplication $\T^\times=\T\setminus\{0_{\T}\}\approx\R$
is an honest commutative group with the unit $``1_{\T}"=0$.
Furthermore, the addition and multiplication satisfy to the
distribution law $``a(b+c)"=``ab+ac"$, $a,b,c\in\T$.
These operations may be viewed as a result of the so-called
{\em dequantization} of the classical arithmetic operations
that underlies the {\em patchworking} construction,
see \cite{Li} and \cite{Vi}.

These tropical operations allow one to define
tropical Laurent polynomials. Namely, a tropical Laurent polynomial
is a function $f:\R^n\to\R$,
$$f(x)=``\sum\limits_j a_j x^j"=\max\limits_j (a_j+ j x),$$
where $j x$ denotes the scalar product, $x\in (\T^\times)^n\approx\R^n$,
$j\in \Z^n$ and only finitely
many coefficients $a_j\in\T$ are non-zero (i.e. not $-\infty$).

Affine-linear functions with integer slopes (for brevity we call them
simply {\em affine functions}) form an important
subcollection of all Laurent polynomials. Namely, these are
such functions $f:\R^n\to\R$ that both $f$ and $``\frac{1_T}{f}"=-f$
are tropical Laurent polynomials. 

We equip $\T^n\approx [-\infty,\infty)^n$ with the Euclidean topology.
Let $U\subset\T^n$ be an open set.
\begin{defn} A continuous function $f:U\to\T$ is called {\em regular}
if its restriction to $U\cap\R^n$ coincides with a restriction of
some tropical Laurent polynomial to $U\cap\R^n$.
\end{defn}
We denote the sheaf of regular functions on $\T^n$ with $\OO$
(or sometimes $\OO_{\T^n}$ to avoid confusion).
Any subset $X\subset\T^n$ gets an induced
regular sheaf $\OO_X$ by restriction. For our purposes we
restrict our attention only to the case when
$X$ is a polyhedral complex, i.e. when $X$ is
the closure of a union of convex polyhedra (possibly unbounded)
in $\R^n$ such that the intersection of any number of such
polyhedra is their common face. We say that $X$ is an $k$-dimensional
polyhedral complex if it is obtained from a union of $k$-dimensional
polyhedra. These polyhedra are called the {\em facets} of $X$.

Let $V\subset X$ be an open set and $f\in\OO_{X}(V)$
be a regular function in $V$. A point $x\in V$ is called a ``zero point"
of $f$ if the restriction of $``\frac{1_{\T}}{f}"=-f$ to $W\subset V$
is not regular for any open neighborhood $W\ni x$.
Note that it may happen that $x$ is a ``zero point" for $\phi:U\to\T$,
but not for $\phi|_{X\cap U}$.
It is easy to see that if $X$ is a $k$-dimensional polyhedral complex
then the ``zero locus" $Z_{f}$ of $f$ is a $(k-1)$-dimensional
polyhedral subcomplex.

To each facet of $Z_{f}$ we may associate a natural number,
called its {\em weight} (or degree). To do this we choose a
``zero point" $x$ inside such a facet.
We say that $x$ is a ``simple zero" for $f$
if for any local decomposition into a sum (i.e. the tropical product)
of regular function $f=``gh"=g+h$ on $V$ near $x$ we have either $g$ or $h$ affine
(i.e. without a ``zero"). We say that the weight is $l$
if $f$ can be locally decomposed into a tropical product of $l$
functions with a simple zero at $x$.

A regular function $f$ allows us to make the following modification on
its domain $V\subset X\subset \T^n$. Consider the graph
$$\Gamma_f\subset V\times\T\subset\T^{n+1}.$$
It is easy to see that the ``zero locus" $$\bar\Gamma_f\subset V\times\T$$
of the (regular) function $``y+f(x)"$ (defined on $V\times\T$),
where $x$ is the coordinate on $V$
and $y$ is the coordinate on $\T$, coincides with the union of $\Gamma_f$
and the undergraph $$U\Gamma_{f,Z}=\{(x,y)\in V\times \T\ |\ x\in Z_f, y\le f(x)\}.$$
Furthermore, the weight of a facet $F\subset\bar\Gamma_f$ is 1
if $F\in\Gamma_f$ (recall that as $V$ is an unweighted polyhedral
complex all the weights of its facets are equal to one) and is
the weight of the corresponding
facet of $Z_f$ if $F\in U\Gamma_{f,Z}$.
We view $\bar\Gamma_f$ as a ``tropical closure" of the set-theoretical
graph $\Gamma_f$. Note that we have a map $\bar\Gamma_f\to V$.
We set $\tilde{V}=\bar\Gamma_f$ to be the result of the {\em tropical modification}
$\mu_f:\tilde{V}\to V$ along the regular function $f$.
The locus $Z_f$ is called the {\em center} of tropical modification.

The weights of the facets of $\tilde{V}$ supplies us with some inconvenience
as they should be incorporated in the definition of the regular
sheaf $\OO_{\tilde{V}}$ on $\tilde{V}$.
Namely, the affine functions defined by $\OO_{\tilde{V}}$ on
a facet of weight $w$ should contain
the group of functions that come as restrictions to this
facet of the affine functions on $\T^{n+1}$
as a subgroup of index $w$.

Sometimes
one can get rid of the weights of $\tilde{V}$ by
a reparameterization with the help of a map $\bar{V}\to\tilde{V}$
that is given by locally linear maps in the corresponding charts.
Indeed, the restriction of $\mu_g:\bar{V}\to\tilde{V}$ to a facet is locally
given by a linear function between two $k$-dimensional affine-linear spaces
defined over $\Z$. If its determinant equals to $w$ then the push-forward
of $\OO_{\bar{V}}$ supplies an extension of $\OO_{\tilde{V}}$ required by the weights.
Note however that if $w>1$ then the converse map is not defined over $\Z$ and
thus is not given by elements of $\OO_{\tilde{V}}$.

Tropical modifications give the basic equivalence relation in Tropical Geometry.
It can be shown that if we start from $\T^k$ and do a number of tropical modifications
on it then the result is a $k$-dimensional polyhedral complex $Y\subset \T^n$
that satisfies to the following {\em balancing property} (cf.
Property 3.3 in \cite{Mi-ICM} where balancing is restated in an equivalent way).
\begin{property}\label{prop-b}
Let $E\subset Y\cap\R^N$ be a $(k-1)$-dimensional face and $F_1,\dots,F_l$ be
the facets of $Y$ adjacent to $F$ whose weights are $w_1,\dots,w_l$.
Let $L\subset\R^N$ be a $(N-k)$-dimensional affine-linear space
with an integer slope and such that it intersects $E$.
For a generic (real) vector $v\in\R^N$ the intersection
$F_j\cap (L+v)$ is either empty or a single point.
Let $\Lambda_{F_j}\subset \Z^N$ be the integer vectors parallel to $F_j$
and $\Lambda_L\subset \Z^N$ be the integer vectors parallel to $L$.
Set $\lambda_j$ to be the product of $w_j$ and the index
of the subgroup $\Lambda_{F_j} +\Lambda_L\subset\Z^N$.
We say that $Y\subset\T^n$ is {\em balanced} if for any choice of $E$, $L$
and a small generic $v$ the sum
$$\iota_L=\sum\limits_{j\ |\ F_j\cap (L+v)\neq\emptyset}\lambda_j$$
is independent of $v$.
We say that $Y$ is {\em simply balanced} if in addition for every $j$
we can find $L$ and $v$ so that $F_j\cap (L+v)\neq\emptyset$,
$\iota_L=1$ and for every small $v$ there exists an affine
hyperplane $H_v\subset L$ such that the intersection $Y\cap (L+v)$ sits
entirely on one side of $H_v+v$ in $L+v$ while the intersection
$Y\cap (H_v+v)$ is a point.
\end{property}

\begin{defn}[cf. \cite{Mi-book},\cite{Mi-ICM}]
A topological space $X$ enhanced with a sheaf of tropical functions $\OO_X$
is called a (smooth) tropical variety of dimension $k$ if for every $x\in X$
there exist an open set $U\ni x$ and an open set $V$ in a simply balanced
polyhedral complex $Y\subset \T^N$ such that the restrictions
$\OO_X|_U$ and $\OO_Y|_V$ are isomorphic.
\end{defn}

Tropical varieties are considered up to the equivalence generated
by tropical modifications.
It can be shown that a smooth tropical variety of dimension $k$
can be locally obtained from $\T^k$ by a sequence of tropical
modifications centered at smooth tropical varieties of dimension $(k-1)$.
This follows from the following proposition.

\begin{prop}
Any $k$-dimensional simply balanced polyhedral complex $X\subset\R^n$
can be obtained from $\T^k$ by a sequence of consecutive tropical
modifications whose centers are simply balanced $(k-1)$-dimensional
polyhedral complexes.
\end{prop}

\begin{proof}
We prove this proposition inductively by $n$. Without the loss of genericity
we may assume that $X$ is a fan, i.e. each convex polyhedron of $X$
is a cone centered at the origin.

The base of the induction,
when $n=k$, is trivial. If $n>k$ let us take a $(n-k)$-dimensional
affine-linear subspace $L\subset\R^n$ given by Property \ref{prop-b}.
Choose a linear projection $$\lambda:\R^n\to\R^{n-1}$$ defined over $\Z$
and such that $\ker(\lambda)$ is a line contained in $L$.

The image $\lambda(X)\subset\R^{n-1}$ is a $k$-dimensional polyhedral
complex since $L$ is transversal to some facets of $X$.
We claim that $$\lambda|_X:X\to\lambda(X)$$
is a tropical modification once we identify $\R^n$ and $\R^{n-1}\times\R$.
The center of this modification is the locus
$$Z_f=\{x\in \R^{n-1}\ |\ \dim(\lambda^{-1}(x)\cap X)>0\}.$$
Here we use the dimension in the usual topological sense.
Note that the $(k-1)$-dimensional complex $Z_f\subset\R^{n-1}$
is simply balanced, existence of the needed $(n-k)$-dimensional
affine-linear spaces follows from the fact that $X\subset\R^n$ is
simply balanced.

To justify our claim we note that near any point $x\in Z_f$
the subcomplex $Y\subset X$ obtained as the
(Euclidean topology) closure of $X\setminus\lambda^{-1}(Z_f)$
is a (set-theoretical) graph of a convex function. This, once again,
follows from the fact that $X\subset\R^n$ is
simply balanced, this time applied to the points in the facets on $X\setminus Y$.
Thus it gives a regular tropical function $f$ and it
remains only to show that the the weight of any facet of $E\subset Z_f$ is 1.
But this follows, in turn, from the balancing condition at $\lambda^{-1}(E)\cap Y$.
\end{proof}

\section{Tropical curves and their moduli spaces}
The definition of tropical variety is especially easy in dimension 1.
Tropical modifications take a graph into a graph (with arbitrary
valence of its vertices) and the tropical structure carried
by the sheaf $\OO_X$ amounts to a complete metric on the complement
of the set of 1-valent vertices of the graph $X$ (cf. \cite{Mi-book}, \cite{MZ},
\cite{GM}). Thus, each 1-valent vertex of a tropical curve $X$ is adjacent
to an edge of infinite length.

A tropical modification allows one to contract such an edge or to attach it
at any point of $X$ other than a 1-valent vertex. If we have a finite
collection of marked points on $X$ then by passing to an equivalent
model if needed we may assume that the set of marked points coincides
with the set of 1-valent vertices. (Of course, if $X$ is a tree then we
have to have at least two marked points to make such assumption.)

The {\em genus} of a tropical curve $X$ is $\dim H_1(X)$.
Let $\MM_{g,n}$ be the set of all tropical curves $X$ of genus $g$ with
$n$ distinct marked points. Fixing a combinatorial type of a graph $\Gamma$
with $n$ marked leaves defines a subset $U_{\Gamma}\subset\MM_{g,n}$ consisting
of marked tropical curves with this combinatorics. A length of any non-leaf
edge of $\Gamma$ defines a real-valued function on $U_\Gamma$. Such functions
are called {\em edge-length functions}. To avoid difficulties caused by
self-automorphisms of $X$ from now on we restrict our attention to the
case $g=0$.

\begin{defn}
The {\em combinatorial type} of a tropical curve $X$ is its
equivalence class up to homeomorphisms respecting
the markings.
\end{defn}

Combinatorial types partite the set $\M_{0,n}$ into disjoint subsets.
The edge-length functions define the structure of the polyhedral cone
$\R_{\ge 0}^M$ in each of those subsets (as the lengths have to be positive).
The number $M$ here is the number of the bounded (non-leaf) edges in $X$.
By the Euler characteristic reasoning it is equal to $n-3$ if $X$ is (1- and)
3-valent, it is smaller if $X$ has vertices of higher valence.

Furthermore, any face of the polyhedral cone $\R_{\ge 0}^M$ coincides
with the cone corresponding to another combinatorial type, the one
where we contract some of the edges of $X$ to points. This gives the
adjacency (fan-like) structure on $\MM_{0,n}$, so $\MM_{0,n}$
is a (non-compact) polyhedral complex. In particular, it is a topological space.

\begin{theorem}\label{thm1}
The set $\MM_{0,n}$ for $n\ge 3$ admits the structure of an $(n-3)$-dimensional
tropical variety such that the edge-length functions are regular within
each combinatorial type.
Furthermore, the space $\MM_{0,n}$ can be tropically embedded in $\R^N$ for some $N$
(i.e. $\MM_{0,n}$ can be presented as a simply balanced complex).
\end{theorem}
\begin{proof}
This theorem is trivial for $n=3$ as $\MM_{0,3}$ is a point.
Otherwise, any two disjoint ordered pairs of marked points can be
used to define a global regular function on $\MM_{0,n}$ with values
in $\R=\T^\times$. Namely, each such ordered pair defines
the oriented path on the tropical curve $X$ connecting the corresponding
marked points. These paths can be embedded.

Since the two pairs of marked points are disjoint the intersection
of the two corresponding paths has to have finite length. We take
this length with the positive sign if the orientations agree and
with the negative sign otherwise. This defines a function on $\MM_{0,n}$.
We call such functions the {\em double ratio} functions.

Take all possible disjoint pairs
of marked points and use them as coordinates for our embedding
$$\iota: \MM_{0,n}\to \R^N,$$
where $N$ is the number of all possible decompositions of $n$ into two
disjoint pairs. The theorem now follows from the following two lemmas.
\end{proof}

Note that, strictly speaking, each coordinate in $\R^N$ depends
not only on the choice of two disjoint pairs of marked points but
also on the order of points in each pair. However, changing the order
in one of the pairs only reverses the sign of the double ratio. Taking
an extra coordinate for such a change of order would be redundant.
Indeed, for any balanced complex $Y\subset\R^N$ and any affine-linear function
$\lambda:\R^N\to\R$ with an integer slope the graph of $\lambda$ is a balanced
complex in $\R^{N+1}$ isomorphic to the initial complex $Y$.

\begin{lem}
The map $\iota$ is a topological embedding.
\end{lem}
\begin{proof}
First, let us prove that $\iota$ is an embedding. The combinatorial type
of $X$ is determined by the set of the coordinates
that do not vanish on $X$. Indeed, any non-leaf edge $E$ of the tree $X$
separates the leaves (i.e. the set of markings) into two classes corresponding
to the components of $X\setminus E$. Let us take a coordinate in $\R^n$
that corresponds to four marking points (union of the two disjoint pairs)
such that two of these points belong to one class and two to the other class.
We call such a coordinate an {\em $E$-compatible coordinate}.
Note that an $E$-compatible coordinate
vanishes on $X$ if and only if the pairs
of markings defined by the coordinate agree with the pairs
defined by the classes.

This observation suffices to reconstruct the combinatorial
type of $X$. Furthermore, the length of $E$ equals to the
minimal non-zero absolute value of the $E$-compatible coordinates.
This implies that $\iota$ is an embedding.
\end{proof}

\begin{lem}
The image $\iota(\MM_{0,n})$
is a simply balanced complex in $\R^N$.
\end{lem}
\begin{proof}
This is a condition
on codimension 1 faces of $\MM_{0,n}$. First we shall check it for the case $n=4$.
There are three ways to split the four marking points into two disjoint
pairs. Accordingly, there are three combinatorial types of 3-valent
trees with three marked leaves. Thus our space $\MM_{0,4}$ is homeomorphic
to the {\em tripod}, or the ``interior" of the letter $Y$, see Figure
\ref{m04}. Each ray of this tripod correspond to a combinatorial type
of a 3-valent tree with 4 leaves while the vertex correspond to the
4-valent tree.
\begin{figure}[h]
\centerline{\psfig{figure=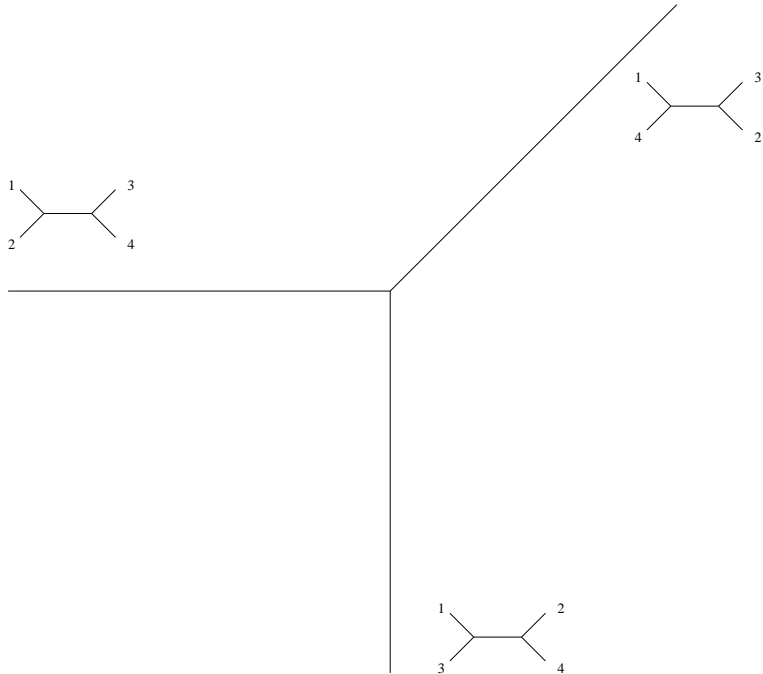}}
\caption{\label{m04} The tropical moduli space $\MM_{0,4}$ and its points
on the corresponding edges.}
\end{figure}

Up to the sign we have the total of three double ratios for $n=4$.
Let us e.g. take those defined by the following ordered pairs:
$\{(12),(34)\}$, $\{(13),(24)\}$ and $\{(14),(23)\}$
Each is vanishing on the corresponding ray of the tripod.
Let us parameterize each ray of the tripod by its only edge-length $t\ge 0$
and compute the corresponding map to $\R^3$.

We have the following embeddings on the three rays
$$t\mapsto (0,t,t),\ t\mapsto (t,0,-t),\ t\mapsto (-t,-t,0).$$
The sum of the primitive integer vectors parallel to the resulting directions
is $0$ and thus $\iota(\MM_{0,4})$ is balanced.

In the case $n>4$ the codimension 1 faces of $\MM_{0,n}$
correspond to the combinatorial types of $X$ with a single 4-valent vertex.
Near a point inside of such face $F$ the space
$\MM_{0,n}$ looks like the product of $\MM_{0,4}$ and $\R^{n-4}$.
The factor $\R^{n-4}$ comes from the edge-lengths on $F$ (its
combinatorial type has $n-4$ bounded edges) while the factor
$\MM_{0,4}$ comes from perturbations of the 4-valent vertex
(which result in a new bounded edge in one of the three possible
combinatorial types of the result).

We have a well-defined map from the union $U$ of the $F$-adjacent facets to $F$
by contracting the new edge to a point. Note that the edge-length functions
exhibit $F$ as the positive quadrant in $\R^{n-4}$.
Furthermore, in the combinatorial type of $F$ we may choose
4 leaves such that contracting all other leaves will take place
outside of the 4-valent vertex (see Figure \ref{contraction4}).
This contraction defines a map $U\to\MM_{0,4}$.
\begin{figure}[h]
\centerline{\psfig{figure=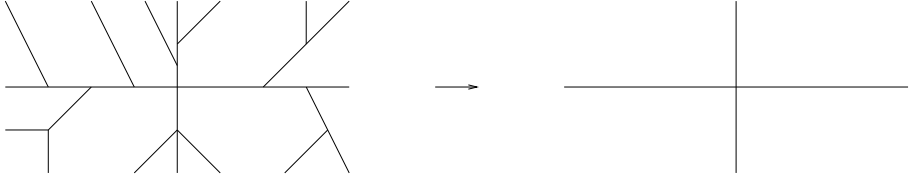}}
\caption{\label{contraction4} One of the possible contractions of a tree
with a 4-valent vertex to the tree corresponding to the origin $O\in\MM_{0,4}$.}
\end{figure}

The lemma now follows from the observation that the resulting decomposition
into $\MM_{0,4}\times\R^{n-4}$ agrees with the double ratio functions.
Indeed, note that the complement of the 4-valent vertex for a curve in
the combinatorial type $F$ is composed of four components.
If the double ratio is such that its four markings are in one-to-one
correspondence with these components then at $U$ it coincides with
sum of the pull-back of the corresponding
double ratio in $F$ with the pull-back of the corresponding
double ratio in $\MM_{0,4}$. If one of the four components is lacking
a marking from the double ratio $\rho$ then $\rho|_U$ coincides
with the corresponding pull-back from $F$.
\end{proof}

\begin{rmk}
The functions $Z_{x_i,x_j}$ from \cite{Mi-ICM} do not define
regular functions on $\MM_{0,n}$, contrary to what is written in \cite{Mi-ICM}.
These functions were a result of an erroneous simplification of the double ratio
functions. But these functions cannot be regular as they are always positive
and Proposition 5.12 of \cite{Mi-ICM} is not correct.
Even the projectivization of the embedding is not a balanced complex already
for $\MM_{0,5}$.
One should use the (non-simplified) double ratios instead.
\end{rmk}

Clearly, the space $\MM_{0,n}$ is non compact. However it is easy
to compactify it by allowing the lengths of bounded edges to assume infinite
values. Let $\overline{\MM}_{0,n}$ be the space of connected trees
with $n$ (marked) leaves such that each edge of this tree
is assigned a length $0<l\le+\infty$ so that each leaf has length
necessarily equal to $+\infty$.

\begin{cor}
The space $\overline{\MM}_{0,n}$ is a smooth compact tropical variety.
\end{cor}

To verify that $\overline{\MM}_{0,n}$ is smooth near a point $x$
at the boundary
$$\dd\overline{\MM}_{0,n}=\overline{\MM}_{0,n}\setminus \MM_{0,n}$$
we need to examine those double ratios that are equal to $\pm\infty$
at $x$. There we use only those signs that result in $-\infty$ do that
the map takes values in $\T^N$.

\begin{rmk}
Note that the compactification $\overline{\MM}_{0,n}\supset \MM_{0,n}$
corresponds to the Deligne-Mumford compactification in the complex case
as under the 1-parametric family collapse of a Riemannian surface
to a tropical curve the tropical length of an edge corresponds to the
rate of growth of the complex modulus of the holomorphic
annulus collapsing to that edge.

Furthermore, similarly to the complex story
the infinite edges decompose a tropical curve into components
(where the non-leaf edges are finite). Any tropical map from
an infinite edge which is bounded would have to be constant
and thus the image would have to split as a union of several
tropical curves in the target.
Such decompositions were used by Gathmann and Markwig in their
deduction of the tropical WDVV equation in $\R^2$, see \cite{GM}.
\end{rmk}

\section{Tropical $\psi$-classes}
Note that we do have the forgetting maps
$$\ft_j:\overline{\MM}_{0,n+1}\to\overline{\MM}_{0,n}$$
for $j=1,\dots,n+1$ by contracting the leaf with the $j$-marking.
This map is sometimes called the {\em universal curve}.
Each marking $k\neq j$ defines a section $\sigma_k$ of $\ft_j$.
The conormal bundle to $\sigma_k$ defines the $\psi_k$-class
in complex geometry (to avoid ambiguity we take $j=n+1$).
This notion can be adapted to our tropical setup.

Recall that so far our choice of tropical models in their equivalence class
was such that the leaves of the tropical curves were
in 1-1 correspondence with the markings.
For this choice we have the images $\sigma_k(\overline{\MM}_{0,n})$ contained
in the boundary part of $\overline{\MM}_{0,n+1}$.
This presentation is compatible with the point of view
when we think about line bundles in tropical geometry
to be given by $H^1(X,\OO^\times)$. Here $X$ is the base of the bundle
and $\OO^\times$ is the sheaf of ``non-vanishing" tropical regular functions.
Such functions are given in the charts to $\R^N$ by affine-linear functions
with integer slopes, see \cite{MZ}.
(Recall that $\T^\times=\R$ is an honest group with respect to tropical multiplication,
i.e. the classical addition.)

However, the following alternative construction allows one
to obtain the $\psi$-classes more geometrically (as we'll illustrate
in an example in the next section).
This approach is based on contracting the leaves marked by number $k$.

The canonical class of a tropical curve is supported at its vertices,
namely we take each vertex with the multiplicity equal to its valence minus 2,
cf. \cite{MZ}.
Furthermore, the cotangent bundle near a 3-valent vertex point
can be viewed as a neighborhood of the origin for the line given
by the tropical polynomial $``x+y+1_{\T}"$ in $\R^2$, so
the $+1$ self-intersection of the line gives the required multiplicity
for the canonical class at any 3-valent vertex.
Thus we can use the intersections with the corresponding codimension 1 faces
in $\MM_{0,n}$ to define the $\psi$-classes there. In other words,
tropical $\psi$-classes will be supported on the $(n-4)$-dimensional
faces in $\MM_{0,n}$.

Namely, for a $\psi_k$-class we have
to collect those codimension 1 faces in $\MM_{0,n}$
whose only 4-valent vertex is adjacent
to the leaf marked by $k$. After a contraction of this leaf
we get a 3-valent vertex, thus the multiplicity of every face
in a $\psi$-divisor is 1. We arrive to the following definition.

\begin{defn}
The tropical $\psi_k$-divisor $\Psi_k\subset\MM_{0,n}$
is the union of those $(n-4)$-dimensional faces that correspond
to tropical curves with a 4-valent vertex adjacent to the leaf marked by $k$,
$k=1,\dots,n$. Each such face is taken with the multiplicity 1.
\end{defn}

\begin{prop}
The subcomplex $\Psi_k$ is a divisor, i.e. satisfies the balancing
condition.
\end{prop}
\begin{proof}
Recall that the balancing condition is a condition at $(n-5)$-dimensional
faces. In $\MM_{0,n}$ there are two types of such faces, one
corresponding to tropical curves with two 4-valent vertices
and one corresponding to a tropical curve with a 5-valent vertex.

Near the faces of the first type the moduli space
$\MM_{0,n}$ is locally a product of
two copies of $\MM_{0,4}$ and $\R^{n-5}$. The $\Psi$-divisor is
a product of $\R^{n-5}$, one copy of $\MM_{0,4}$ and the central (3-valent)
point in the other copy of $\MM_{0,4}$ (this is the point corresponding
to the 4-valent vertex adjacent to the leaf marked by $k$). Thus
the balancing condition holds trivially in this case.

Near the faces of the second type the moduli space
$\MM_{0,n}$ is locally a product of $\MM_{0,5}$ and $\R^{n-5}$.
As in the proof of Theorem \ref{thm1} each double ratio decomposes
to the sum of the corresponding double ration in $\MM_{0,5}$
(perhaps trivial if two of the markings for the double ratio
correspond to the same edge adjacent to the 5-valent vertex)
and an affine-linear function in $\R^{n-5}$.
Thus it suffices to check only the balancing condition for the
$\Psi$-divisors in $\MM_{0,5}$. This example is considered in details
in the next section. The balancing condition there follows from
Proposition \ref{bal5}.
\end{proof}

Conjecturally, the tropical $\Psi$-divisors are limits of
some natural representatives of the divisors
for the complex $\psi$-classes under the collapse of the
complex moduli space onto the corresponding tropical moduli space
$\MM_{0,n}$. Note that our choice for the tropical $\Psi$-divisor
is not contained in the boundary $\dd\overline{\MM}_{0,n}\subset
\overline{\MM}_{0,n}$ (cf. the calculus of the complex
boundary classes in \cite{Keel}), but comes as a closure
of a divisor in $\MM_{0,n}$.

\section{The space $\overline{\MM}_{0,5}$}
We have already described the moduli space $\MM_{0,4}$ as the
tripod of Figure \ref{m04}. It has only one $0$-dimensional face $O\in\MM_{0,4}$.
This point (considered as a divisor) coincides with the divisors
$\Psi_1=\Psi_2=\Psi_3=\Psi_4$. The description of $\MM_{0,5}$
is somewhat more interesting.

\begin{figure}[h]
\centerline{\psfig{figure=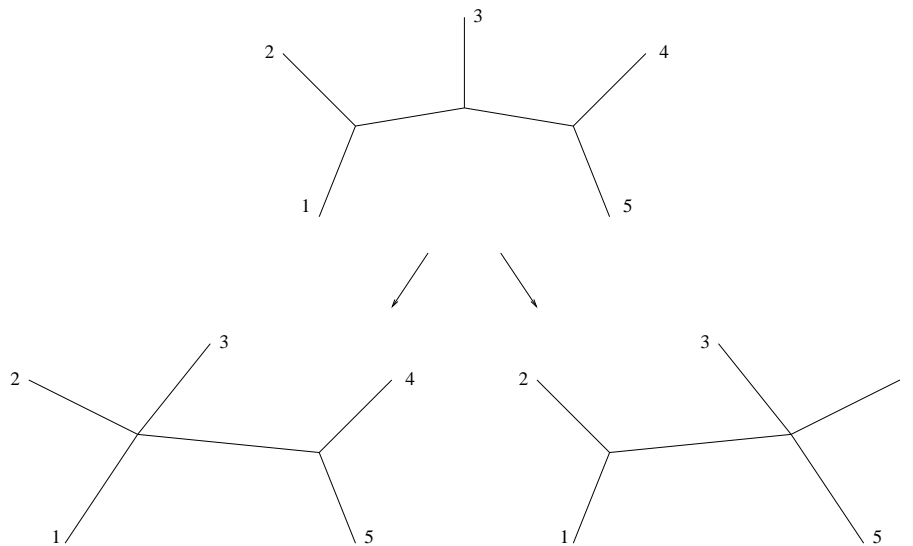}}
\caption{\label{tree5} Adjunction of combinatorial types corresponding
to the quadrant connecting the rays $(45)$ and $(12)$.}
\end{figure}

There are 15 combinatorial types of 3-valent trees with 5 marked
leaves. If we forget about the markings there is only one homeomorphism
class for such a curve (see Figure \ref{tree5}). To get the number
of non-isomorphic markings we take the number all possible reordering of vertices
(equal to $5!=120$) and divide by $2^3=8$ as there is an 8-fold symmetry
of reordering. Indeed there is one symmetry interchanging the left two
leaves, one interchanging the right two leaves and the central symmetry
around the central leave of the 3-valent tree on top of Figure \ref{tree5}.
\begin{figure}[h]
\centerline{\psfig{figure=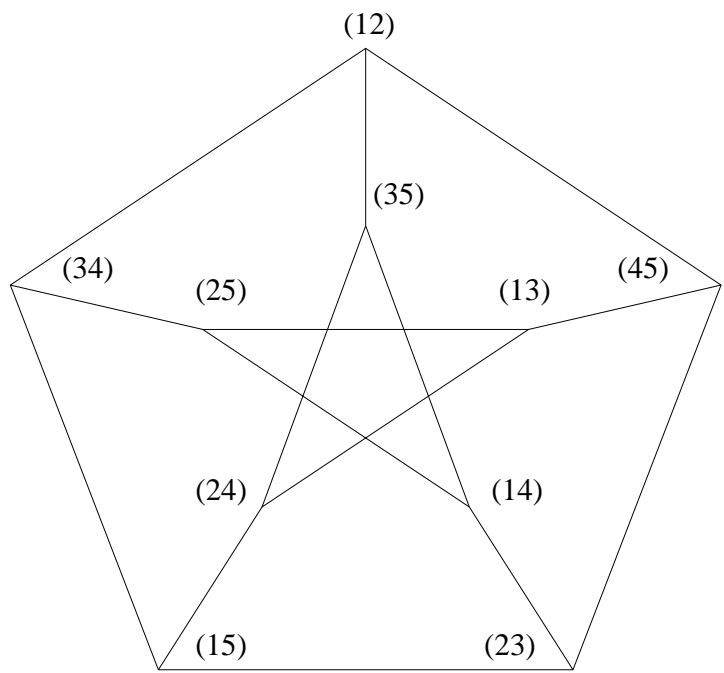}}
\caption{\label{m05} The link of the origin in $\MM_{0,5}$.}
\end{figure}

Thus the space $\MM_{0,5}$ is a union of 15 quadrants $\R^2_{\ge 0}$.
These quadrants are attached along the rays which correspond to
the combinatorial types of curves with one 4-valent vertex.
Such curves also have one 3-valent vertex which is adjacent to
two leaves and the only bounded edge of the curve, see the bottom
of Figure \ref{tree5}. Such combinatorial
types are determined by the markings of the two leaves
emanating from the 3-valent vertex. Thus we have a total
of $\begin{pmatrix} 5 \\ 2 \end{pmatrix}=10$ of such rays.

The two boundary edges of the quadrant correspond to contractions
of the bounded edges of the combinatorial type as shown on Figure \ref{tree5}.
The global picture of adjacency of quadrants and rays is shown
on Figure \ref{m05} where the reader may recognize the well-known
{\em Petersen graph}, cf. the related tropical Grassmannian picture
in \cite{SS}. Vertices of this graph correspond to the rays of
$\MM_{0,5}$ while the edges correspond to the quadrants. Thus
the whole picture may be interpreted as the link of the only
vertex $O\in\MM_{0,5}$ (the point $O$ corresponds to the tree with a
5-valent vertex adjacent to all the leaves).

To locate the $\Psi_k$-divisor we recall that the $k$th leaf has
to be adjacent to a 4-valent vertex if it appears in $\Psi_k$.
This means that $\Psi_k$ consists of 6 rays that are marked by pairs
not containing $k$.

\begin{prop}\label{bal5}
The subcomplex $\Psi_k\subset\MM_{0,5}$ is a divisor.
\end{prop}
\begin{proof}
Since the whole $\MM_{0,5}$ is $S_5$-symmetric
it suffices to check the balancing condition only for $\Psi_1$.
The embedding $\MM_{0,5}\subset\R^N$ is given by the double ratios,
so it suffices to check that for each double ratio function
the sum of its gradients on the six rays of $\Psi_1$ vanishes.

If the double ratio is determined by two pairs disjoint from the marking 1,
e.g. by $\{(23),(45)\}$ then its restriction onto the six rays of $\Psi_1$
is the same as its restriction to the three rays $\MM_{0,4}$ taken twice
and thus balanced. Namely its gradient is 1 on the rays $(24)$ and $(35)$;
$-1$ on the rays $(25)$ and $(34)$; and 0 on the rays $(23)$ and $(45)$.

If the four markings of the double ratio contain the marking 1
then thanks to the symmetry we may assume that the double ratio
is given by $\{(12),(34)\}$. It vanishes on the rays $(34)$, $(35)$,
$(45)$ and $(25)$; it has gradient $+1$ on the ray $(24)$ and the gradient
$-1$ on the ray $(23)$. Once again, the balancing condition holds.
\end{proof}

\begin{figure}[h]
\centerline{\psfig{figure=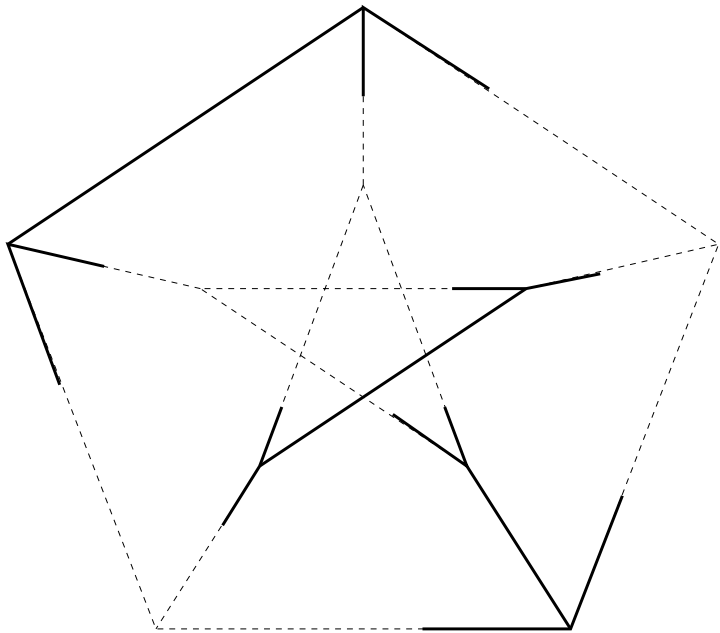}}
\caption{\label{unicurve} The three fibers and four sections
of the universal curve $\ft_5:\MM_{0,5}\to\MM_{0,4}$.}
\end{figure}

As our final example of the paper we would like to describe explicitly
the universal curve $$\ft_5:\MM_{0,5}\to\MM_{0,4}.$$
This is presented on Figure \ref{unicurve}.
Once again, we interpret the Peterson graph as the link $L$ of
the vertex $O\in\MM_{0,5}$. Similarly, the link of the origin
in $\MM_{0,4}$ consists of three points.
Thus $L$ is the union of the fibers of $\ft_5$ (away from a
neighborhood of infinity) over these three points and four copies of
a neighborhood of the origin in $\MM_{0,4}$ corresponding to the four sections
$\sigma_1$, $\sigma_2$, $\sigma_3$ and $\sigma_4$ of the universal curve.
Figure \ref{unicurve} depicts the fibers in $L$ with solid lines and the
sections with dashed lines.

\begin{ack}
I am thankful to Valery Alexeev and Kristin Shaw for discussions
related to geometry of tropical moduli spaces. My research is supported
in part by NSERC.
\end{ack}


\begin{thebibliography}{99}

\bibitem{GM} Gathmann, A., Markwig, H.,
Kontsevich's formula and the WDVV equations in tropical geometry,
http://arxiv.org/abs/math.AG/0509628.

\bibitem{Keel} Keel, S., Intersection theory of moduli space of
stable $N$-pointed curves of genus zero, Transactions of the AMS
{\bf 330} (1992), 545--574.

\bibitem{Li}  Litvinov, G. L.,
The Maslov dequantization, idempotent and tropical mathematics: a very brief introduction.
In {\em Idempotent mathematics and mathematical physics},
Contemp. Math., 377, Amer. Math. Soc., Providence, RI, 2005, 1--17.

\bibitem{Mi-ICM} Mikhalkin, G., {Tropical Geometry and its application},
to appear in the Proceedings on the ICM-2006, Madrid;
http://arxiv.org/abs/math/0601041.

\bibitem{Mi-book} Mikhalkin, G., {Tropical Geometry},
book in preparation.

\bibitem{MZ} Mikhalkin, G., Zharkov, I.,
Tropical curves, their Jacobians and Theta functions,
http://arxiv.org/abs/math/0612267.

\bibitem{SS} Speyer, D., Sturmfels, B.,
The tropical Grassmannian.  {\em Adv. Geom.}  4  (2004),  no. 3, 389--411.

\bibitem{Vi}  Viro, O. Ya., Dequantization of real algebraic geometry on logarithmic paper.
In {\em European Congress of Mathematics, Vol. I (Barcelona, 2000),}
Progr. Math., 201, Birkh\"auser, Basel, 2001, 135--146.

\end{thebibliography}
\end{document}